\documentclass[a4paper]{amsart}
\usepackage{epsfig}
\usepackage[utf8]{inputenc} 
\usepackage[all]{xy}
\usepackage{amssymb}
\usepackage{amsmath}
\newtheorem{theorem}{Theorem}

{\it}{\rm}

\newtheorem*{bem}{Remark}{\it}{}

\numberwithin{equation}{section}

\newcommand{\Z}{{\mathbb Z}}

\newcommand{\bH}{{\mathbb H}}

\newcommand{\reteil}{\operatorname*{Re}}

\begin{document}

\title[]
{Restricted sums of four integral squares}  
\author[R. Schulze-Pillot]{Rainer Schulze-Pillot} 
 \begin{abstract}
We give a simple quaternionic proof of a recent result of Goldmakher
and Pollack on restricted sums of four integral squares.
 \end{abstract}
 \maketitle
In \cite{gopo} the authors prove the following result:
\begin{theorem}[Goldmakher, Pollack]
 Let $n \equiv
 T \bmod 2$ be integers. Then $n$ has a representation
 $n=\sum_{\nu=0}^3a_\nu^2$ as a sum of four integer squares with
 $\sum_{\nu=0}^3a_\nu=T$ if and only if $4n-T^2$ is a sum of three
 integral squares. 
\end{theorem}
We give here a different proof using a very simple computation in the
ring of integral quaternions.
\begin{proof}
Let $1,i,j,k$ denote the usual basis of the Hamilton quaternions
$\bH$, for $\alpha=a_0+a_1i+a_2j+a_3k \in \bH$ write $\reteil(\alpha):=a_0,
N(\alpha):=\sum_{\nu=0}^3a_\nu^2, \varphi(\alpha):=\sum_{\nu=0}^3
a_\nu$. We put $R=\Z+\Z i+\Z j+\Z k$ and define $f:R\to R$ by
$f(\alpha)=\alpha(1-i-j-k)$. We have $\reteil(f(\alpha))=\varphi(\alpha), N(f(\alpha))=4N(\alpha)$,
and since $N(\alpha)\equiv \varphi(\alpha)\bmod 2$ for $\alpha \in R$, we
have 
$f(R)\subseteq X:=\{\beta \in R \mid N(\beta)\equiv 4\reteil(\beta)
\bmod 8\}$.

Since $(1-i-j-k)^{-1}=(1+i+j+k)/4$ we see that
$f(R)=\{\beta \in R\mid \beta(1+i+j+k)/4 \in R\}=\{\beta=b_0+b_1i+b_2j+b_3k\in R\ \mid \sum_{\nu\ne \mu}b_\nu
\equiv b_\mu \bmod 4 \text{ for } 0\le \mu \le 3\}$. 

For
$\beta'=b_0'+b_1'i+b_2'j+b_3'k \in  X$ one easily sees that the
$b_\nu'$ are all congruent modulo $2$ and if they are even, the number
of $b_\nu'\equiv 2 \bmod 4$ is $0,2$ or $4$. In the latter case, the
congruence condition $\sum_{\nu\ne \mu}b_\nu'
\equiv b_\mu' \bmod 4 \text{ for } 0\le \mu \le 3$ is satisfied. If
the $b_\nu'$ are all odd, it is satisfied if either one or three of
them are congruent to $-1$ modulo $4$.
In both cases we find, changing a sign if necessary,  a
$\beta=b_0+b_1i+b_2j+b_3k \in f(R)$ with $b_0=b_0', b_\nu\in
\{b_\nu',-b_\nu'\}$ for all $\nu$. 

Obviously, for $\alpha=a_0+a_1i+a_2j+a_3k \in R$ with
$\sum_{\nu=0}^3a_\nu^2=n, \sum_{\nu=0}^3
a_\nu=T$ and $f(\alpha)=:b_0+b_1i+b_2j+b_3k$ we have
$\sum_{\nu=1}^3 b_\nu^2=4n-T^2$, which establishes the ``only if''
part of the
the assertion.

Conversely, let $n\equiv T \bmod 2$ 
be such that $4n-T^2=\sum_{\nu=1}^3(b_\nu')^2$ is a sum of three
squares.
For $\beta'=T+b_1'i+b_2'j+b_3'k$  we have $N(\beta)\equiv 4
\reteil(\beta')\bmod 8$, hence $\beta'\in X$. 
By what we have shown about elements of $X$ there exists $\beta=T+b_1i+b_2j+b_3k \in f(R)$
with $N(\beta)=N(\beta')=4n$, and $\alpha =a_0+a_1i+a_2j+a_3k \in R$ with
$\beta=f(\alpha)$ is as required. 

\end{proof}
\begin{bem}
  Our proof shows more precisely that $r_4(n,T)=\vert \{{\bf a}\in
  \Z^4\mid \sum a_\nu^2=n, \sum a_\nu=T\}\vert$  equals
  $r_3(4n-T^2)=\vert\{{\bf b}\in \Z^3\mid \sum b_\nu^2=4n-T^2\}\vert$
  if $n$ and $T$ are even and $\frac{r_3(4n-T^2)}{2}$ if $n$ and $T$
  are odd. This can also be proved with the help of identities for
  Jacobi theta series.
\end{bem}

Rainer Schulze-Pillot,
Fachrichtung 
Mathematik,
Universit\"at des Saarlandes\\
Postfach 151150, 66041 Saarbr\"ucken, Germany\\
email: schulzep@math.uni-sb.de

\end{document}